# An $L^1$-type estimate for Riesz potentials

Armin Schikorra[*1], Daniel Spector[†2], and Jean Van Schaftingen[‡3]

[1]Department of Mathematics, University of Basel, Basel, Switzerland
[2]Department of Applied Mathematics, National Chiao Tung University, Hsinchu, Taiwan
[3]Institut de Recherche en Mathématique et en Physique, Université catholique de Louvain, Louvain-la-Neuve, Belgium

**Abstract**

In this paper we establish new $L^1$-type estimates for the classical Riesz potentials of order $\alpha \in (0, N)$:

$$\|I_\alpha u\|_{L^{N/(N-\alpha)}(\mathbb{R}^N)} \leq C \|Ru\|_{L^1(\mathbb{R}^N;\mathbb{R}^N)}.$$

This sharpens the result of Stein and Weiss on the mapping properties of Riesz potentials on the real Hardy space $\mathcal{H}^1(\mathbb{R}^N)$ and provides a new family of $L^1$-Sobolev inequalities for the Riesz fractional gradient.

## 1 Introduction and Main Results

Let $N \geq 2$ and define the Riesz potential $I_\alpha$ of order $\alpha \in (0, N)$ by its action on a measurable function $u$ via the convolution

$$I_\alpha u(x) \equiv (I_\alpha * u)(x) := \frac{1}{\gamma(\alpha)} \int_{\mathbb{R}^N} \frac{u(y)}{|x-y|^{N-\alpha}} \, dy,$$

whenever it is well-defined. Here, $\gamma(\alpha) = \pi^{N/2} 2^\alpha \Gamma(\frac{\alpha}{2})/\Gamma(\frac{N-\alpha}{2})$ is a normalization constant [24, p. 117] that ensures that the Riesz potentials satisfy the semigroup property

$$I_{\alpha+\beta} u = I_\alpha I_\beta u, \qquad \text{for } \alpha, \beta > 0, \text{ such that } \alpha + \beta < N,$$

for $u$ in a suitable class of functions.

The study of the mapping properties of $I_\alpha$ on $L^p(\mathbb{R}^N)$ was initiated by Sobolev, who proved the following *fundamental theorem about integrals of the potential type* in 1938 [21, p. 50].

[*]armin.schikorra@unibas.ch, A.S. supported by SNF
[†]dspector@math.nctu.edu.tw, D.S. supported by MOST 103-2115-M-009-016-MY2
[‡]Jean.VanSchaftingen@uclouvain.be

**Theorem 1.1** (Sobolev)**.** *Let $0 < \alpha < N$ and $1 < p < N/\alpha$. Then there exists a constant $C = C(p, \alpha, N) > 0$ such that*

$$\|I_\alpha u\|_{L^{Np/(N-\alpha p)}(\mathbb{R}^N)} \leq C\|u\|_{L^p(\mathbb{R}^N)} \tag{1.1}$$

*for all $u \in L^p(\mathbb{R}^N)$.*

In particular we see that Sobolev's result concerns $L^p$ estimates for Riesz potentials when $1 < p < \frac{N}{\alpha}$ and strictly excludes the case $p = 1$. Indeed, it is well-known that no such inequality as (1.1) can hold in this regime - one may consider, for example (cf. [24, p. 119]), an approximation of the identity (for explicit construction one can see Section 3 of this paper). Then the right-hand-side of the inequality in Theorem 1.1 stays bounded while pointwise

$$I_\alpha \rho_\varepsilon(x) \to I_\alpha(x) = \frac{1}{\gamma(\alpha)} \frac{1}{|x|^{N-\alpha}},$$

which does not belong to the Lebesgue space $L^{\frac{N}{N-\alpha}}(\mathbb{R}^N)$, and Fatou's lemma gives the desired contradiction.

It is then natural to ask if there is a substitute for the inequality (1.1). One possibility was given by Stein and Weiss [25, p. 31], where they demonstrated that if one replaces $L^p(\mathbb{R}^N)$ with the real Hardy space

$$\mathcal{H}^p(\mathbb{R}^N) := \{u \in L^p(\mathbb{R}^N) : Ru \in L^p(\mathbb{R}^N; \mathbb{R}^N)\}$$

(where $Ru := DI_1 u$ is the vector-valued Riesz transform), one can extend the validity of Theorem 1.1 to the regime $p = 1$. For $p \in (1, \infty)$, $\mathcal{H}^p(\mathbb{R}^N) = L^p(\mathbb{R}^N)$, but for $p = 1$ the Hardy space $\mathcal{H}^1(\mathbb{R}^N)$ is strictly contained in $L^1(\mathbb{R}^N)$. Their result implies the following theorem.

**Theorem 1.2** (Stein-Weiss)**.** *Let $0 < \alpha < N$ and $1 \leq p < N/\alpha$. Then there exists a constant $C = C(p, \alpha, N) > 0$ such that*

$$\|I_\alpha u\|_{L^{Np/(N-\alpha p)}(\mathbb{R}^N)} \leq C \left(\|u\|_{L^p(\mathbb{R}^N)} + \|Ru\|_{L^p(\mathbb{R}^N;\mathbb{R}^N)}\right)$$

*for all $u \in \mathcal{H}^p(\mathbb{R}^N)$.*

Actually the approach to Sobolev inequalities due to Gagliardo [10, p. 120] and Nirenberg [17, p. 128] gives another replacement to Theorem 1.1 for $1 \leq \alpha < N$. Indeed, written in the language of potentials, one sees that the results [10,17] assert the existence of a constant $C > 0$ such that

$$\|I_1 u\|_{L^{N/(N-1)}(\mathbb{R}^N)} \leq C\|Ru\|_{L^1(\mathbb{R}^N;\mathbb{R}^N)},$$

for all $u \in C_c^\infty(\mathbb{R}^N)$ such that $Ru \in L^1(\mathbb{R}^N; \mathbb{R}^N)$. Therefore, if $1 \leq \alpha < N$, the preceding inequality and Theorem 1.1 applied to $I_\alpha u = I_{\alpha-1} I_1 u$ allows us to deduce that

$$\|I_\alpha u\|_{L^{N/(N-\alpha)}(\mathbb{R}^N)} \leq C\|I_1 u\|_{L^{N/(N-1)}(\mathbb{R}^N)}$$
$$\leq C\|Ru\|_{L^1(\mathbb{R}^N;\mathbb{R}^N)},$$

for all $u \in C_c^\infty(\mathbb{R}^N)$ such that $Ru \in L^1(\mathbb{R}^N; \mathbb{R}^N)$.

The main result of this paper is the following theorem demonstrating that this $L^1$-type estimate holds for the Riesz potential of any order $\alpha \in (0, N)$.

**Theorem A.** Let $N \geq 2$ and $0 < \alpha < N$. Then there exists a constant $C = C(\alpha, N) > 0$ such that

$$\|I_\alpha u\|_{L^{N/(N-\alpha)}(\mathbb{R}^N)} \leq C \|Ru\|_{L^1(\mathbb{R}^N;\mathbb{R}^N)}$$

for all $u \in C_c^\infty(\mathbb{R}^N)$ such that $Ru \in L^1(\mathbb{R}^N;\mathbb{R}^N)$.

**Remark 1.3.** Theorem A is false when $N = 1$, see Counterexample 3.2 in Section 3.

Our motivation for such an inequality can be found in the study of certain fractional partial differential equations introduced in [20], where existence results are demonstrated for a continuous spectrum of such equations parameterized by the Riesz fractional gradient

$$D^\alpha u := D I_{1-\alpha} u,$$

for $0 < \alpha < 1$. With this notation, an alternative formulation of Theorem A is the following.

**Theorem A′.** Let $N \geq 2$ and $0 < \alpha < 1$. Then there exists a constant $C = C(\alpha, N) > 0$ such that

$$\|u\|_{L^{N/(N-\alpha)}(\mathbb{R}^N)} \leq C \|D^\alpha u\|_{L^1(\mathbb{R}^N;\mathbb{R}^N)} \tag{1.2}$$

for all $u \in C_c^\infty(\mathbb{R}^N)$.

Theorem A′ is a natural analogy to the Sobolev inequalities known for the fractional Laplacian when $p > 1$ and integer order derivatives for $p \geq 1$, though one might have guessed such a theorem from several additional factors. Firstly, related results for Besov spaces with the same degree of fractional differentiability have long been known in the literature (see e.g. [6, Lemma D.2; 8, Theorem 1.4; 13, Theorem 4; 22, Theorem 2; 32, Theorem 8.3]). A second factor suggesting such an inequality is the observation that the asymptotics of the constant in Theorem 1.1 are $O(1/(p-1))$ as $p \to 1$, which agrees with the asymptotics of the operator norm of the vector-valued Riesz transform $R : L^p(\mathbb{R}^N) \to L^p(\mathbb{R}^N; \mathbb{R}^N)$. Finally, there is the more recent work of the second author and R. Garg [11, 12] which shows that the logarithmic potential $I_N u$, defined for $u \in C_c^\infty(\mathbb{R}^N)$ by

$$I_N u(x) = \frac{1}{|S^{N-1}|} \int_{\mathbb{R}^N} \log\left(\frac{1}{|x-y|}\right) u(y)\, dy,$$

has, for any $u$ with[1] $\int_{\mathbb{R}^N} u = 0$, the representation

$$I_N u(x) = \frac{1}{|S^{N-1}|} \int_{\mathbb{R}^N} \frac{x-y}{|x-y|} \cdot Ru(y)\, dy.$$

Therefore, when $\alpha = N$ one has the corresponding estimate

$$\|I_N u\|_{L^\infty(\mathbb{R}^N)} \leq C \|Ru\|_{L^1(\mathbb{R}^N;\mathbb{R}^N)}.$$

---

[1] For this class of functions, this is equivalent to asking $Ru \in L^1(\mathbb{R}^N; \mathbb{R}^N)$.

Our proof of Theorem A is quite direct, and relies only on the boundedness of the Riesz transform and of the classical maximal function operator on $L^p$ for $1 < p < +\infty$. We do not rely upon any Sobolev type embedding nor any multiplier theorem that goes beyond the Riesz transform. Here the crucial observation is that the vector-valued Riesz transform is curl-free, i.e.

$$\frac{\partial R_j u}{\partial x_i} = \frac{\partial R_i u}{\partial x_j}$$

for all $i, j \in \{1, \ldots, N\}$. In fact, an interesting point to note is that the same proof shows that one has

$$\|I_\alpha F\|_{L^{N/(N-\alpha)}(\mathbb{R}^N;\mathbb{R}^N)} \leq C\|F\|_{L^1(\mathbb{R}^N;\mathbb{R}^N)} \tag{1.3}$$

for vector fields $F \in L^1(\mathbb{R}^N; \mathbb{R}^N)$ that satisfy either $\operatorname{curl} F = 0$ or $\operatorname{div} F = 0$, the pair of which is reminiscent of the conditions for inclusion in the real Hardy space [25].

The remainder of the paper is organized as follows. In Section 2 we give proofs of the main results. In Section 3 we discuss several more intricate questions in greater detail, including connections of our result with more technical results from the literature, an open question in regard to a sharp result in the scale of Lorentz spaces, and the details of the counterexample mentioned in Remark 1.3.

## 2   Proofs of the Main Results

We now prove Theorem A. In the course of the proof, we will use $C$ to designate a constant that may depend upon $\alpha$ and $N$, though the constant may change from line to line.

*Proof of Theorem A.* Let $u \in C_c^\infty(\mathbb{R}^N)$ be such that $Ru \in L^1(\mathbb{R}^N; \mathbb{R}^N)$. We claim it suffices to show that, for $j \in \{1, \ldots, N\}$, one has the existence of a uniform constant $C = C(\alpha, N) > 0$ such that

$$\left|\int_{\mathbb{R}^N} R_j u I_\alpha \varphi\right| \leq C\|Ru\|_{L^1(\mathbb{R}^N;\mathbb{R}^N)} \|\varphi\|_{L^{\frac{N}{\alpha}}(\mathbb{R}^N)} \tag{2.1}$$

for every $\varphi \in C_c^\infty(\mathbb{R}^N)$.

Indeed, utilizing the identity $v = \sum_{i=1}^N -R_i^2 v$, the boundedness of $R_i : L^{N/(N-\alpha)}(\mathbb{R}^N) \to L^{N/(N-\alpha)}(\mathbb{R}^N)$, and duality we have

$$\|I_\alpha u\|_{L^{N/(N-\alpha)}(\mathbb{R}^N)} \leq C \sum_{i=1}^N \|R_i I_\alpha u\|_{L^{N/(N-\alpha)}(\mathbb{R}^N)}$$

$$= C \sum_{i=1}^N \sup_{\|\psi_i\|_{L^{N/\alpha}} \leq 1} \int_{\mathbb{R}^N} R_i I_\alpha u\, \psi_i.$$

Now, Fubini's theorem implies that

$$\int_{\mathbb{R}^N} R_i I_\alpha u\, \psi_i = \int_{\mathbb{R}^N} R_i u\, I_\alpha \psi_i.$$

Therefore we can estimate

$$\|I_\alpha u\|_{L^{N/(N-\alpha)}(\mathbb{R}^N)} \leq C \sum_{i=1}^{N} \sup_{\|\psi_i\|_{L^{N/\alpha}} \leq 1} \|Ru\|_{L^1(\mathbb{R}^N;\mathbb{R}^N)} \|\psi_i\|_{L^{\frac{N}{\alpha}}(\mathbb{R}^N)}$$

$$= C \sum_{i=1}^{N} \sup_{\|\psi_i\|_{L^{N/\alpha}} \leq 1} \|Ru\|_{L^1(\mathbb{R}^N;\mathbb{R}^N)} \|\psi_i\|_{L^{\frac{N}{\alpha}}(\mathbb{R}^N)}$$

$$\leq C \|Ru\|_{L^1(\mathbb{R}^N;\mathbb{R}^N)},$$

which is the thesis.

We therefore proceed to prove inequality (2.1). We follow the strategy of [28]. Without loss of generality, we take $j = 1$ and write $x = (x', x_N)$. We now introduce a family of mollifiers: we take $\rho \in C_c^\infty(\mathbb{R}^N)$ such that $\operatorname{supp} \rho \subset B(0,1)$, and

$$\int_{\mathbb{R}^N} \rho = 1.$$

Then we define $\rho_\varepsilon(x) = \rho(x/\varepsilon)/\varepsilon^N$ and $\varphi_\varepsilon(x) = (\varphi * \rho_\varepsilon)(x)$ so that

$$\int_{\mathbb{R}^{N-1}} R_1 u(x', x_N) \, I_\alpha \varphi(x', x_N) \, dx' = \int_{\mathbb{R}^{N-1}} R_1 u(x', x_N) \, [I_\alpha \varphi(x', x_N) - I_\alpha \varphi_\varepsilon(x', x_N)] \, dx'$$

$$+ \int_{\mathbb{R}^{N-1}} R_1 u(x', x_N) \, I_\alpha \varphi_\varepsilon(x', x_N) \, dx'$$

$$=: I(\varepsilon) + II(\varepsilon).$$

For the first term we begin with the bound

$$I(\varepsilon) \leq \|R_1 u(\cdot, x_N)\|_{L^1(\mathbb{R}^{N-1})} \|I_\alpha \varphi(\cdot, x_N) - I_\alpha \varphi_\varepsilon(\cdot, x_N)\|_{L^\infty(\mathbb{R}^{N-1})}.$$

Now, the fundamental theorem of calculus implies that

$$I_\alpha \varphi_\varepsilon(x) - I_\alpha \varphi_\delta(x) = \int_\delta^\varepsilon \frac{\partial}{\partial r} \rho_r * I_\alpha \varphi(x) \, dr$$

$$= \int_\delta^\varepsilon \int_{\mathbb{R}^N} \sigma_r(x-y) \, I_\alpha \varphi(y) \, dy dr,$$

where

$$\sigma_r(z) := \frac{\partial \rho_r}{\partial r}(z) = \frac{1}{r^N} \left[ -\nabla \rho\left(\frac{z}{r}\right) \cdot \frac{z}{r^2} - \frac{N}{r} \rho\left(\frac{z}{r}\right) \right].$$

Thus sending $\delta \to 0$, Lebesgue's dominated convergence theorem implies

$$I_\alpha \varphi_\varepsilon(x) - I_\alpha \varphi(x) = \int_0^\varepsilon \int_{\mathbb{R}^N} \sigma_r(x-y) \, I_\alpha \varphi(y) \, dy dr.$$

As before, Fubini's theorem yields the identity

$$I_\alpha \varphi_\varepsilon(x) - I_\alpha \varphi(x) = \int_0^\varepsilon \int_{\mathbb{R}^N} I_\alpha \sigma_r(x-y) \varphi(y) \, dy dr.$$

Next, we claim that one has the pointwise inequality

$$|I_\alpha \sigma_r(z)| \leq \frac{C}{(r+|z|)^{N-\alpha+1}}. \tag{2.2}$$

We distinguish two cases: $|z| \leq 2r$ and $|z| > 2r$. When $|z| \leq 2r$, one has

$$|I_\alpha \sigma_r(z)| = \frac{C}{r^N} \left| \int_{B(0,r)} \frac{\nabla \rho\left(\frac{z}{r}\right) \cdot \frac{z}{r^2} + \frac{N}{r} \rho\left(\frac{z}{r}\right)}{|z-y|^{N-\alpha}} \, dy \right|$$

$$\leq \frac{C}{r^{N+1}} \int_{B(0,r)} \frac{1}{|z-y|^{N-\alpha}} \, dy$$

$$\leq \frac{C}{r^{N-\alpha+1}}.$$

Then $|z| \leq 2r$ implies $1/r^{N-\alpha+1} \leq \frac{3^{N-\alpha+1}}{(|z|+r)^{N-\alpha+1}}$, which allows us to deduce the inequality (2.2) in this regime. Next, when $|z| > 2r$, we have

$$I_\alpha \sigma_r(z) = \frac{C}{r^N} \int_{B(0,r)} \frac{\operatorname{div}\left(\rho\left(\frac{y}{r}\right) \frac{y}{r}\right)}{|z-y|^{N-\alpha}} \, dy$$

$$= \frac{1}{r^N} \int_{B(0,r)} -\frac{\rho\left(\frac{y}{r}\right) \frac{y}{r}}{|z-y|^{N-\alpha+1}} \cdot \frac{y-z}{|y-z|} \, dy,$$

which upon the change of variables $w = y/r$ yields the bound

$$|I_\alpha \sigma_r(z)| \leq \int_{B(0,1)} \frac{C}{|z-rw|^{N-\alpha+1}} \, dw$$

$$= \frac{1}{|z|^{N-\alpha+1}} \int_{B(0,1)} \frac{C}{\left|\frac{z}{|z|} - \frac{r}{|z|} w\right|^{N-\alpha+1}} \, dw.$$

Finally, the assumption that we are in the regime $|z| > 2r$ implies both that the last integral is bounded and in a similar manner to before that $1/|z|^{N-\alpha+1} \leq \frac{C}{(|z|+r)^{N-\alpha+1}}$, thus proving (2.2).

Therefore, we can estimate

$$\|I_\alpha \varphi(\cdot, x_N) - I_\alpha \varphi_\varepsilon(\cdot, x_N)\|_{L^\infty(\mathbb{R}^{N-1})} \leq C \sup_{x' \in \mathbb{R}^{N-1}} \int_0^\varepsilon \int_{\mathbb{R}^N} \frac{|\varphi(y)|}{(r+|x-y|)^{N-\alpha+1}} \, dy \, dr.$$

By the Hölder inequality on $\mathbb{R}^{N-1}$, we deduce that

$$\|I_\alpha \varphi(\cdot, x_N) - I_\alpha \varphi_\varepsilon(\cdot, x_N)\|_{L^\infty(\mathbb{R}^{N-1})}$$

$$\leq C \sup_{x' \in \mathbb{R}^{N-1}} \int_0^\varepsilon \int_{\mathbb{R}} \left( \int_{\mathbb{R}^{N-1}} |\varphi(y', y_N)|^{\frac{N}{\alpha}} \, dy' \right)^{\frac{\alpha}{N}}$$

$$\left( \int_{\mathbb{R}^{N-1}} \frac{1}{(r+\sqrt{|x_N-y_N|^2 + |x'-y'|^2})^{N+\frac{N}{N-\alpha}}} \, dy' \right)^{1-\frac{\alpha}{N}} dy_N \, dr.$$

If we set

$$\Phi(y_N) = \left( \int_{\mathbb{R}^{N-1}} |\varphi(y', y_N)|^{\frac{N}{\alpha}} \, dy' \right)^{\frac{\alpha}{N}},$$

and can establish the estimate

$$\sup_{x' \in \mathbb{R}^{N-1}} \left( \int_{\mathbb{R}^{N-1}} \frac{1}{(r+\sqrt{|x_N-y_N|^2 + |x'-y'|^2})^{N+\frac{N}{N-\alpha}}} \, dy' \right)^{1-\frac{\alpha}{N}}$$

$$\leq \frac{C}{(r+|x_N-y_N|)^{2-\frac{\alpha}{N}}},$$

then we would have the bound

$$\|I_\alpha\varphi(\cdot,x_N) - I_\alpha\varphi_\varepsilon(\cdot,x_N)\|_{L^\infty(\mathbb{R}^{N-1})}$$
$$\leq C \int_0^\varepsilon \int_\mathbb{R} \frac{\Phi(y_N)}{(r+|x_N-y_N|)^{1-\frac{\alpha}{N}}} \, dy_N \, dr.$$

However, let us observe that

$$r + |x_N - y_N| + |x' - y'| \leq C\left(r + \sqrt{|x_N-y_N|^2 + |x'-y'|^2}\right),$$

and so

$$\int_{\mathbb{R}^{N-1}} \frac{1}{(r+\sqrt{|x_N-y_N|^2+|x'-y'|^2})^{N+\frac{N}{N-\alpha}}} \, dy'$$
$$\leq \int_{\mathbb{R}^{N-1}} \frac{C}{(r+|x_N-y_N|+|x'-y'|)^{N+\frac{N}{N-\alpha}}} \, dy'$$
$$= \frac{C}{(r+|x_N-y_N|)^{N+\frac{N}{N-\alpha}}} \int_{\mathbb{R}^{N-1}} \frac{1}{(1+\frac{|x'-y'|}{r+|x_N-y_N|})^{N+\frac{N}{N-\alpha}}} \, dy'.$$

Integrating in spherical coordinates with center at $x'$, we find

$$\int_{\mathbb{R}^{N-1}} \frac{1}{(1+\frac{|x'-y'|}{r+|x_N-y_N|})^{N+\frac{N}{N-\alpha}}} \, dy'$$
$$= C\int_0^\infty \frac{t^{N-2}}{(1+\frac{t}{r+|x_N-y_N|})^{N+\frac{N}{N-\alpha}}} \, dt$$
$$= C(r+|x_N-y_N|)^{N-1} \int_0^\infty \frac{(t')^{N-2}}{(1+t')^{N+\frac{N}{N-\alpha}}} \, dt',$$

from which the result follows.

Finally, considering the integrand on dyadic annuli we have

$$\int_0^\varepsilon \int_\mathbb{R} \frac{\Phi(y_N)}{(r+|x_N-y_N|)^{2-\frac{\alpha}{N}}} \, dy_N \, dr$$
$$= \int_0^\varepsilon \sum_{n\in\mathbb{Z}} \int_{2^n r<|x_N-y_N|<2^{n+1}r} \frac{\Phi(y_N)}{(r+|x_N-y_N|)^{2-\frac{\alpha}{N}}} \, dy_N \, dr$$
$$\leq \int_0^\varepsilon \sum_n \frac{2^{n+1}r}{(1+2^n)^{2-\frac{\alpha}{N}}} \frac{1}{r^{2-\frac{\alpha}{N}}} \fint_{B(x_N,2^{n+1}r)} \Phi(y_N) \, dy_N \, dr$$
$$\leq \sum_n \frac{2^{n+1}}{(1+2^n)^{2-\frac{\alpha}{N}}} \left(\int_0^\varepsilon \frac{1}{r^{1-\frac{\alpha}{N}}} \, dr\right) \mathcal{M}\Phi(x_N),$$

where $\mathcal{M}\Phi : \mathbb{R} \to \mathbb{R}$ is the Hardy–Littlewood maximal function of $\Phi : \mathbb{R} \to \mathbb{R}$. We have thus

$$I(\varepsilon) \leq C\|R_1 u(\cdot,x_N)\|_{L^1(\mathbb{R}^{N-1})} \varepsilon^{\frac{\alpha}{N}} \mathcal{M}\Phi(x_N). \tag{2.3}$$

Now for $II(\varepsilon)$ we apply the fundamental theorem of calculus to write

$$II(\varepsilon) = -\int_{\mathbb{R}^{N-1}} \int_{x_N}^\infty \frac{\partial}{\partial x_N} R_1 u(x',t) \, I_\alpha\varphi_\varepsilon(x',x_N) \, dt \, dx'.$$

We use the fact that the vector-valued Riesz transform is curl-free, i.e.

$$\frac{\partial R_j u}{\partial x_i} = \frac{\partial R_i u}{\partial x_j}$$

for all $i, j \in \{1, \ldots, N\}$ and Fubini's theorem to deduce that

$$-\int_{\mathbb{R}^{N-1}} \int_{x_N}^{\infty} \frac{\partial}{\partial x_N} R_1 u(x', t) \, I_\alpha \varphi_\varepsilon(x', x_N) dt dx'$$

$$= -\int_{\mathbb{R}^{N-1}} \int_{x_N}^{\infty} \frac{\partial}{\partial x_1} R_N u(x', t) \, I_\alpha \varphi_\varepsilon(x', x_N) dt dx'$$

$$= \int_{x_N}^{\infty} -\int_{\mathbb{R}^{N-1}} \frac{\partial}{\partial x_1} R_N u(x', t) \, I_\alpha \varphi_\varepsilon(x', x_N) dx' dt.$$

The important point now is that $N \neq 1$, allowing us to integrate by parts and obtain

$$\int_{x_N}^{\infty} -\int_{\mathbb{R}^{N-1}} \frac{\partial}{\partial x_1} R_N u(x', t) \, I_\alpha \varphi_\varepsilon(x', x_N) dx' dt$$

$$= \int_{x_N}^{\infty} \int_{\mathbb{R}^{N-1}} R_N u(x', t) \, \frac{\partial}{\partial x_1} I_\alpha \varphi_\varepsilon(x', x_N) dx' dt.$$

Thus,

$$II(\varepsilon) \leq \|R_N u\|_{L^1(\mathbb{R}^N)} \sup_{x' \in \mathbb{R}^{N-1}} \left| \frac{\partial I_\alpha \varphi_\varepsilon}{\partial x_1}(x', x_N) \right|.$$

In a similar manner to the first case, we see that

$$\frac{\partial I_\alpha \varphi_\varepsilon(x', x_N)}{\partial x_1} = \int_{\mathbb{R}^N} \frac{\partial \rho_\varepsilon(y)}{\partial x_1} I_\alpha \varphi(x - y) \, dy$$

$$= \int_{\mathbb{R}^N} I_\alpha \frac{\partial \rho_\varepsilon}{\partial x_1}(y) \varphi(x - y) \, dy,$$

where we again have the pointwise estimate

$$\left| I_\alpha \frac{\partial \rho_\varepsilon}{\partial x_1}(y) \right| \leq \frac{C}{(\varepsilon + |y|)^{N-\alpha+1}}.$$

Therefore, Hölder's inequality in $\mathbb{R}^{N-1}$ with an analogous estimate to the preceding yields the bound

$$\sup_{x' \in \mathbb{R}^{N-1}} \left| \frac{\partial I_\alpha \varphi_\varepsilon}{\partial x_1}(x', x_N) \right| \leq C \int_{\mathbb{R}} \frac{\Phi(x_N - y_N)}{(\varepsilon + |y_N|)^{2-\frac{\alpha}{N}}} \, dy_N.$$

Finally, the dyadic splitting can again be employed to enable one to conclude

$$\sup_{x' \in \mathbb{R}^{N-1}} \left| \frac{\partial I_\alpha \varphi_\varepsilon}{\partial x_1}(x', x_N) \right| \leq C \frac{\mathcal{M}\Phi(x_N)}{\varepsilon^{1-\frac{\alpha}{N}}},$$

so that

$$II(\varepsilon) \leq C \|R_N u\|_{L^1(\mathbb{R}^N)} \frac{\mathcal{M}\Phi(x_N)}{\varepsilon^{1-\frac{\alpha}{N}}}. \tag{2.4}$$

Choosing $\varepsilon = \|R_N u\|_{L^1(\mathbb{R}^N)}/\|R_1 u(\cdot, x_N)\|_{L^1(\mathbb{R}^{N-1})}$, equations (2.3) and (2.4) imply that

$$\int_{\mathbb{R}^{N-1}} R_1 u(x', x_N)\, I_\alpha \varphi(x', x_N) dx' \leq C \|R_1 u(\cdot, x_N)\|_{L^1(\mathbb{R}^{N-1})}^{1-\frac{\alpha}{N}} \|R_N u\|_{L^1(\mathbb{R}^N)}^{\frac{\alpha}{N}} \mathcal{M}\Phi(x_N).$$

We now integrate this estimate with respect to $x_N$ on $\mathbb{R}$ to obtain by the classical Hölder inequality

$$\int_{\mathbb{R}^N} R_1 u\, I_\alpha \varphi \leq C \|R_N u\|_{L^1(\mathbb{R}^N)}^{\frac{\alpha}{N}} \left( \int_{\mathbb{R}} \|R_1 u(\cdot, x_N)\|_{L^1(\mathbb{R}^{N-1})}\, dx_N \right)^{1-\frac{\alpha}{N}}$$
$$\times \left( \int_{\mathbb{R}} \big(\mathcal{M}\Phi(x_N)\big)^{\frac{N}{\alpha}}\, dx_N \right)^{\frac{\alpha}{N}}.$$

By the classical maximal function theorem

$$\int_{\mathbb{R}} \big(\mathcal{M}\Phi(x_N)\big)^{\frac{N}{\alpha}}\, dx_N \leq C \int_{\mathbb{R}} \Phi(x_N)^{\frac{N}{\alpha}}\, dx_N = C \int_{\mathbb{R}^N} |\varphi(x)|^{\frac{N}{\alpha}} dx,$$

which completes the proof of the claim and hence the theorem. $\square$

**Remark 2.1.** Maximal function bounds on the integrals on slices of the type (2.3) and (2.4) were introduced by Chanillo and Van Schaftingen [9].

*Proof of Theorem A'.* Theorem A' can be proven in a similar manner, beginning with an estimate for $u$ in the space $L^{N/(N-\alpha)}(\mathbb{R}^N)$. $\square$

## 3 Connections, Improvements, Counterexamples

### 3.1 Connections to several results in the literature

We have here given a proof of Theorem A (and one can similarly prove Theorem A') using elementary arguments, though there are other possible proofs that could be employed. We mention several here for both historical propriety, and to satisfy the curious reader. In Section 1, we have seen that such a result can be deduced directly from classical and well-known results in the case $\alpha \geq 1$. However, a second method that works for all values of $\alpha \in (0, N)$ can be used if one is willing to accept the embeddings of $W^{1,1}(\mathbb{R}^N)$ and $BV(\mathbb{R}^N)$ into Besov spaces explored in [6, Lemma D.2; 8, Theorem 1.4; 13, Theorem 4; 22, Theorem 2; 32, Theorem 8.3]). One then obtains the result by a combination of these embeddings with the embeddings of Besov spaces into Triebel-Lizorkin spaces, for example if $N \geq 2$ and $\alpha \in (0, 1)$ and denoting

$$\frac{1}{p} = 1 - \frac{\alpha}{N},$$

one has that $1 < p \leq 2$, and so

$$\|v\|_{\dot{F}^{1-\alpha}_{p,2}(\mathbb{R}^N)} \leq C \|v\|_{\dot{B}^{1-\alpha}_{p,p}(\mathbb{R}^N)} \leq C \|Dv\|_{L^1(\mathbb{R}^N;\mathbb{R}^N)},$$

which from a characterization of the space $\dot{F}^{1-\alpha}_{p,2}$ and taking $v = I_1 u$ implies

$$\|I_\alpha u\|_{L^p(\mathbb{R}^N)} \leq C \|Ru\|_{L^1(\mathbb{R}^N;\mathbb{R}^N)},$$

which is the inequality in Theorem A while the inequality in Theorem A' follows in a similar manner taking $v = I_{1-\alpha}u$.

Finally, as in the original proof of the authors, one can argue by duality. This method was pioneered by Bourgain and Brezis in the works [2,3,5], who were interested in constructing bounded solutions to the divergence equation

$$-\operatorname{div} \mathbf{Y} = f$$

in the critical regime $f \in L^N$. The dual result to this is a stronger form of the inequality of Gagliardo and Nirenberg mentioned in the introduction. A simpler proof of this result was given by the third author in [28,31,32], which is the basic idea behind the slicing argument we have utilized. We also mention that when $\alpha > 1/2$ one has a stronger inequality in the spirit of the work of Bourgain and Brezis via the analogous estimates by Bousquet, Mironescu, and Russ [7] in the scale of Triebel-Lizorkin spaces.

### 3.2 Lorentz space improvements

As in the case of embeddings for Sobolev spaces, Theorems A and A' are sharp in the scale of $L^p$ spaces, though can be improved when one considers the finer scale of Lorentz spaces. For instance, in Theorem 1.1 one can replace the $L^{Np/(N-\alpha p)}(\mathbb{R}^N)$ norm on the left hand side with that of the Lorentz space $L^{Np/(N-\alpha p),p}(\mathbb{R}^N)$ (see [18, p. 139]). While $L^{p,p}(\mathbb{R}^N) = L^p(\mathbb{R}^N)$, a smaller second parameter in the Lorentz spaces is more than microscopic improvement. One can easily see this fact by comparing Trudinger's result [27] that $Du \in L^N(\mathbb{R}^N)$ implies $u$ is exponentially integrable (and not in general bounded) with Stein's result [23] that $Du \in L^{N,1}(\mathbb{R}^N)$ implies that $u$ is continuous.

In fact, we can show an estimate in the Lorentz space $L^{N/(N-\alpha),q}(\mathbb{R}^N)$ for any $q > 1$ as follows. Let $q' = N/\varepsilon$ for some $\varepsilon > 0$. Then we utilize inequality (2.1) to find a $C = C(\varepsilon, N) > 0$ such that

$$\Big|\int_{\mathbb{R}^N} R_j u \, I_\epsilon \varphi \Big| \leq C \|Ru\|_{L^1(\mathbb{R}^N;\mathbb{R}^N)} \|\varphi\|_{L^{\frac{N}{\varepsilon}}(\mathbb{R}^N)},$$

which combined with the boundedness of $I_{\alpha-\varepsilon} : L^q(\mathbb{R}^N) \to L^{N/(N-\alpha),q}(\mathbb{R}^N)$ [18, p. 139] implies

$$\|I_\alpha u\|_{L^{N/(N-\alpha),q}(\mathbb{R}^N)} \leq \|I_\varepsilon u\|_{L^q(\mathbb{R}^N)}$$
$$\leq C\|Ru\|_{L^1(\mathbb{R}^N;\mathbb{R}^N)},$$

which gives the desired result.

When $\alpha > 1$, then $I_\alpha u = I_{\alpha-1} I_1 u$, and $\|I_1 u\|_{L^{N/(N-1),1}} \leq C\|Du\|_{L^1}$ [1, 26], which combined with the previously cited convolution estimates of O'Neill [18, p. 139] shows that one can obtain $q = 1$ in this regime. Thus, the critical case here is the endpoint $q = 1$ and $\alpha \in (0,1)$, for which none of the preceding techniques can obviously applied to obtain the estimate. This leads to the following open question concerning a sharper $L^1$-type estimate.

**Open Question 3.1.** Let $N \geq 2$ and suppose $0 < \alpha < 1$. Does there exists a constant $C = C(\alpha, N) > 0$ such that

$$\|I_\alpha u\|_{L^{N/(N-\alpha),1}(\mathbb{R}^N)} \leq C\|Ru\|_{L^1(\mathbb{R}^N;\mathbb{R}^N)}$$

for all $u \in C_c^\infty(\mathbb{R}^N)$ such that $Ru \in L^1(\mathbb{R}^N;\mathbb{R}^N)$?

## 3.3 Counterexamples

We now provide the counterexample mentioned in Remark 1.3, substantiating our claim that Theorem A is false when $N = 1$. Note the similarity to the standard counterexample for the failure of Sobolev's result in $L^1$ discussed in the introduction.

**Counterexample 3.2.** Suppose one had such an inequality as given in Theorem A. Then by density one obtains the inequality for all functions in the real Hardy space $\mathcal{H}^1(\mathbb{R})$. Now, in this setting $Ru = Hu$ is the Hilbert transform, which by our assumption, the identity $H^2 = -I$, and boundedness on $L^{1/(1-\alpha)}(\mathbb{R})$ would imply

$$\|I_\alpha v\|_{L^{1/(1-\alpha)}(\mathbb{R})} \leq C\|HI_\alpha v\|_{L^{1/(1-\alpha)}(\mathbb{R})} \leq C\|v\|_{L^1(\mathbb{R})},$$

for all $v \in \mathcal{H}^1(\mathbb{R})$. Now taking $v_\varepsilon(x) = \rho((x-1)/\varepsilon)/\varepsilon - \rho((x+1)/\varepsilon)/\varepsilon$ with $\rho$ as above ($v_\varepsilon$ is the difference of two translated approximations of the identity), one has $v_\varepsilon \in \mathcal{H}^1(\mathbb{R})$ (as before, for a smooth, compactly supported function a necessary and sufficient condition for this inclusion is that $\int v_\varepsilon = 0$), the right-hand-side stays bounded, and

$$I_\alpha v_\varepsilon(x) \to \frac{1}{\gamma(\alpha)} \left( \frac{1}{|x-1|^{1-\alpha}} - \frac{1}{|x+1|^{1-\alpha}} \right),$$

which fails to be locally $L^{1/(1-\alpha)}$ near $-1$ and $+1$, and so an application of Fatou's lemma gives one the desired contradiction.

The paper [19] contains an example of a one-dimensional failure of an embedding of certain Besov spaces into $BV(\mathbb{R})$. The preceding counterexample combined with the known embeddings for Besov spaces previously discussed gives another such example.

## 4 Acknowledgements

The authors would like to thank Aline Bonami for discussions regarding Riesz potentials and Hardy spaces, Augusto Ponce for critically reviewing the manuscript, as well as the anonymous referee whose comments have improved its presentation. The first author is supported by the SNF and would like to thank D. Spector, NCTU ("Regularity Estimates for Fractional PDE" - Research Grant 104W986) and the CMMSC for their hospitality during the work of this project. The second author is supported by the Taiwan Ministry of Science and Technology under research grant MOST 103-2115-M-009-016-MY2.